\documentclass[twoside,12pt]{article}
\usepackage[latin2]{inputenc}
\usepackage{amssymb}
\usepackage{paralist}
\usepackage{amsmath}
\usepackage{amsfonts}
\usepackage{graphicx}
\usepackage{wrapfig}
\usepackage{tikz}
\usepackage{hyperref}

\usetikzlibrary{arrows}
\usetikzlibrary{decorations.pathmorphing}
\usetikzlibrary{snakes}

\newtheorem{prop}{Proposition}[section]
\newtheorem{thm}{Theorem}[section]

\newtheorem{cor}{Corollary}[prop]

\newtheorem{conj}{Open problem}

\newcommand{\ve}[1]{\mathbf{#1}}
\newcommand{\F}{\mathbb{F}}
\newcommand{\Fx}{\F[x_1,\dots,x_n]}

\newcommand{\tab}{\ \ \ \ }

\date{}

\begin{document}

\vspace{3cm}

\centerline {\Large{\bf Shattering-extremal set systems of small
$VC$-dimension}}

%

\centerline{}

\centerline{\bf {Tam\'as M\'esz\'aros}}

\centerline{}

\centerline{Department of Mathematics, Central European
University}

\centerline{Meszaros\_Tamas@ceu-budapest.edu}


\centerline{}

\centerline{\bf {Lajos R\'onyai\footnote{Research supported in
part by OTKA grants K77476, K77778 and TÁMOP grant
4.2.2.B-10/1-2010-0009.}}}

\centerline{}

\centerline{Computer and Automation Research Institute, Hungarian
Academy of Sciences}

\centerline{Institute of Mathematics, Budapest University of
Technology and Economics}

\centerline{lajos@ilab.sztaki.hu}



\begin{abstract}
We say that a set system $\mathcal{F}\subseteq 2^{[n]}$ shatters a
given set $S\subseteq [n]$ if $2^S=\{F \cap S : F \in
\mathcal{F}\}$. The Sauer inequality states that in general, a set
system $\mathcal{F}$ shatters at least $|\mathcal{F}|$ sets. Here
we concentrate on the case of equality. A set system is called
shattering-extremal if it shatters exactly $|\mathcal{F}|$ sets.
We characterize shattering extremal set systems of
Vapnik-Chervonenkis dimension $1$ in terms of their inclusion
graphs. Also from the perspective of extremality, we relate set
systems of bounded Vapnik-Chervonenkis dimension to their
projections. \\
\end{abstract}



\section{Introduction}
\tab Throughout this paper $n$ will be a positive integer, the set
$\{1,2,\dots,n\}$ will be referred to shortly as $[n]$, the power
set of it as $2^{[n]}$ and the family of subsets of size $k$ as
${[n] \choose k }$.

The central notion of our study is shattering. We say that a set
system $\mathcal{F}\subseteq 2^{[n]}$ shatters a given set
$S\subseteq [n]$ if $2^S=\{F \cap S : F \in \mathcal{F}\}$. The
family of subsets of $[n]$ shattered by $\mathcal{F}$ is denoted
by $Sh(\mathcal{F})$. The following inequality gives a bound on
the size of $Sh(\mathcal{F})$.

\begin{prop} \label{Sauer inequality}
$|Sh(\mathcal{F})|\geq |\mathcal{F}|$.
\end{prop}

The statement was proved by several authors, Aharoni and Holzman
\cite{Aharoni-Holzman}, Pajor \cite{Pajor}, Sauer \cite{Sauer},
Shelah \cite{Shelah}. Often it is referred to as Sauer inequality.
One of the most interesting cases is the case of equality, i.e.
when the set system $\mathcal{F}$ shatters exactly $|\mathcal{F}|$
sets. We call such set systems shattering extremal or s-extremal
for short. Many interesting results have been obtained in
connection with these combinatorial objects, among others in
\cite{Shattering news}, \cite{reverse kleitman}, \cite{defect
sauer}, \cite{Frankl}, and \cite{Furedi}.\\

The Vapnik-Chervonenkis dimension of a set system
$\mathcal{F}\subseteq 2^{[n]}$, denoted by $VC-dim(\mathcal{F})$,
is the maximum cardinality of a set shattered by $\mathcal{F}$.
The general task of giving a good description of s-extremal
systems seems to be too complex at this point. We restrict
therefore our attention to the simplest cases, where the
$VC$-dimension of $\mathcal{F}$ is bounded by some fixed natural
number $t$.

After the introduction, in Section \ref{1 VC dim} we first
investigate s-extremal set systems of $VC$-dimension at most $1$
from a graph theoretical point of view. We give a bijection
between the family of such set systems on the ground set $[n]$ and
trees on $n+1$ vertices. As a consequence one can exactly
determine the number of such s-extremal set systems. In
combinatorics when considering set systems with a given property
it is a common step to first consider families of some special
structure. According to \cite{greco} uniform set systems can't be
s-extremal. As a next possibility set systems from two consecutive
layers turn up. In Section \ref{2 layers} we prove that they are
just special cases of the previous ones. After this in Section
\ref{GR bases 1 VC dim} we switch to an algebraic point of view
and investigate bases the polynomial ideals attached to extremal
set systems. The main result of Section \ref{t VC dim} is a
connection between s-extremal set systems of $VC$-dimension $t$
and their projections. At the end we propose an open problem and
make some concluding remarks.

In the paper we will use the terminology of \cite{DiestelGT} for
graph theoretical notions.

\section{s-extremal set systems of $VC$-dimension at most $1$}\label{1 VC dim}

\tab Let $\mathcal{F}\subseteq 2^{[n]}$ be an s-extremal family.
Let $G_{\mathcal{F}}$ be the labelled Hasse diagram of
$\mathcal{F}$ considered as a graph, i.e. a graph whose vertices
are the elements of $\mathcal{F}$ and there is a directed edge
going from $G$ to $F$, labelled with $j\in [n]$ exactly when
$F=G\cup \{j\}$. $G_{\mathcal{F}}$ will be called the inclusion
graph of $\mathcal{F}$. When representing the elements of
$\mathcal{F}$ by their characteristic vectors, $G_{\mathcal{F}}$
can also be considered as the subgraph in the Hamming graph
$H_n=\{0,1\}^n$ spanned by the elements corresponding to the sets
in $\mathcal{F}$ with edges directed and labelled in a natural
way. Actually for the next proposition we can forget about the
directions of the edges, and consider $G_{\mathcal{F}}$ as an
undirected edge-labelled graph. We further assume that
$VC-dim(\mathcal{F})\leq 1$. Our aim is to characterize these
kinds of s-extremal set systems in terms of their inclusion graph.

\begin{prop}\label{extr char VC-dim=1}
A set system $\mathcal{F}\subseteq 2^{[n]}$ is s-extremal and of
$VC-dimension$ at most $1$ iff $G_{\mathcal{F}}$ is a tree and all
labels on the edges are different.
\end{prop}
\textbf{Proof:} For the 'only if' direction suppose that
$\mathcal{F}$ is s-extremal and $VC-dim(\mathcal{F})\leq 1$.
According to \cite{greco}, Theorem $5$ we know that
$G_{\mathcal{F}}$ must be isometrically embedded into $H_n$ (i.e.
for any two elements $G,F\in \mathcal{F}$ the distance between $G$
and $F$ is the same in $G_{\mathcal{F}}$ and in $H_n$). This means
in particular, that $G_{\mathcal{F}}$ is connected. Next we prove
that all labels on the edges of $G_{\mathcal{F}}$ are different.
Suppose for contradiction, that there are two edges with the same
label. W.l.o.g. we may assume that this label is $1$. Since there
are no two edges going out from a set with the same label, there
are sets $A,B,C,D\in \mathcal{F}$, all different, such that $1\in
A\cap B$, $C=A\backslash \{1\}$ and $D=B\backslash \{1\}$. Since
$A\neq B$, $A\vartriangle B(=(A\backslash B)\cup(B\backslash A))$
is nonempty, so there is an element $1\neq a\in A\vartriangle B$.
W.l.o.g. we may assume that $a\in A\backslash B$. Now

\[\begin{array}{cccc}
  \{1,a\}\cap A=\{1,a\} & \{1,a\}\cap B=\{1\} & \{1,a\}\cap C=\{a\} & \{1,a\}\cap D=\emptyset \\
\end{array}.\]

So $\{1,a\}$ is shattered by $\{A,B,C,D\}$, consequently $\{1,a\}
\in Sh(\mathcal{F})$, contradicting the assumption
$VC-dim(\mathcal{F})\leq 1$.

To finish with this direction note that the fact that all labels
are different implies that $G_{\mathcal{F}}$ is acyclic. Suppose
for contradiction that it is not the case, and $G_{\mathcal{F}}$
contains a cycle. Pick one edge from this cycle and let $a$ be its
label. On the remaining part of the cycle there must be another
edge labelled with $a$, since it connects a set containing $a$
with one not containing $a$. However this is impossible, since all
labels are different. Adding the connectedness of
$G_{\mathcal{F}}$, we obtain that it is actually a tree as wanted.

For the reverse direction suppose that $G_{\mathcal{F}}$ is a tree
and all labels on the edges are different. It is easily seen that
this implies that $G_{\mathcal{F}}$ is isometrically embedded into
$H_n$. (Otherwise a path from a set $A$ to $B$ in
$G_{\mathcal{F}}$ which is not a shortest in $H_n$ would contain
$2$ edges with the same label, corresponding to the addition and
deletion of the same element of $[n]$.)

Now we prove that $VC-dim(\mathcal{F})\leq 1$. Suppose the
contrary, namely that $\mathcal{F}$ shatters a set of size $2$,
e.g. $\{1,2\}$. This means that there are sets $A,B,C,D\in
\mathcal{F}$ such that

\[\begin{array}{cccc}
  \{1,2\}\cap A=\{1,2\} & \{1,2\}\cap B=\{1\} & \{1,2\}\cap C=\{2\} & \{1,2\}\cap D=\emptyset \\
\end{array}.\]

Consider a shortest path in $G_{\mathcal{F}}$ from $A$ to $B$.
Since $2\in A\backslash B$, this shortest path has to contain an
edge labelled with $2$. Repeating this argument for $C$ and $D$
one gets another, different (since on a shortest path between $A$
and $B$ every set contains the element $1$, on the other hand on a
shortest path between $C$ and $D$ none of the sets does) edge with
label $2$, what contradicts the assumption that all labels are
different.

Now we calculate $Sh(\mathcal{F})$. If $i\in [n]$ is not an edge
label, then either all sets $F\in \mathcal{F}$ contain $i$, ore
none of them does. In particular $\{i\}$ is not shattered by
$\mathcal{F}$. Thus $Sh(\mathcal{F})$ consists of $\emptyset$ and
the sets $\{i\}$, where $i$ is an edge label. However all edge
labels are different, so we get that
$|Sh(\mathcal{F})|=|E(G_{\mathcal{F}})|+1=|\mathcal{F}|$ (since
$G_{\mathcal{F}}$ is a tree), i.e. $\mathcal{F}$ is s-extremal.
$\blacksquare$\\[0,1 cm]

Let $\mathcal{F}\subseteq 2^{[n]}$ be an s-extremal family such
that $supp(\mathcal{F})=\cup_{F\in \mathcal{F}} F=[n]$ and
$\cap_{F\in \mathcal{F}} F=\emptyset$. By Proposition \ref{extr
char VC-dim=1} to every s-extremal family of $VC$-dimen\-sion at
most $1$ one can associate a directed edge-labelled tree
$G_{\mathcal{F}}$, all edges having distinct labels. We have seen
that $Sh(\mathcal{F})$ consists of $\emptyset$ and the sets
$\{i\}$, where $i$ is an edge label. On the other hand, since
$\cap_{F\in \mathcal{F}} F=\emptyset$, we also have that
$Sh(\mathcal{F})=\{\emptyset\}\cup \{\{j\}\ |\ j\in
supp(\mathcal{F})=[n]\}$. As a consequence the tree must have $n$
edges and thus $n+1$ vertices, i.e. such an s-extremal family has
$n+1$ elements.

Now conversely suppose that we are given a directed edge-labelled
tree $T$ on $n+1$ vertices with $n$ edges, all having different
labels. This tree at the same time also defines a set system
$\mathcal{F}_T=\{F_v\ |\ v\in T\}$. Take the edges one by one.
When considering an edge with label $s$ going from $u$ to $v$,
then for all vertices $w$ closer to $v$ than to $u$ in the
undirected tree put $s$ into $F_w$. Clearly $T=G_{\mathcal{F}}$
and by the previous proposition $\mathcal{F}$ must be s-extremal.

To illustrate this, consider the following example with $n=5$.

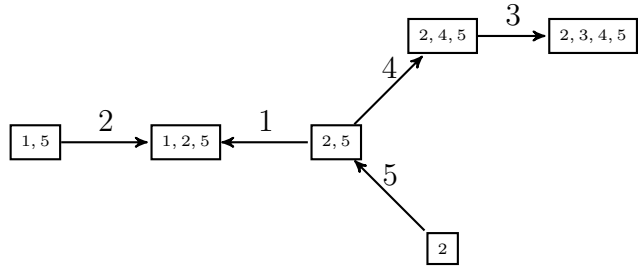
\begin{center}
\begin{tikzpicture}[-,>=stealth',shorten >=1pt,auto,node distance=2cm,
    thick,main node/.style={rectangle,fill=none,draw,font=\sffamily\tiny\bfseries}]

    \node[main node] (1) {$1,5$};
    \node[main node] (2) [right of=1] {$1,2,5$};
    \node[main node] (3) [right of=2] {$2,5$};
    \node[main node] (4) [above right of=3] {$2,4,5$};
    \node[main node] (5) [right of=4] {$2,3,4,5$};
    \node[main node] (6) [below right of=3] {$2$};

    \draw[->] (1) -- (2) node[pos=.5,above] {$2$};
    \draw[<-] (2) -- (3) node[pos=.5,above] {$1$};
    \draw[->] (3) -- (4) node[pos=.5,above] {$4$};
    \draw[->] (4) -- (5) node[pos=.5,above] {$3$};
    \draw[<-] (3) -- (6) node[pos=.5,above] {$5$};

  \end{tikzpicture}
\end{center}

\begin{center}
We have
$\mathcal{F}_T=\{\{1,5\},\{1,2,5\},\{2,5\},\{2,4,5\},\{2,3,4,5\},\{2\}\}$.
\end{center}

This gives a bijection between the set of all s-extremal families
of $VC$-dimension at most $1$ and directed edge-labelled trees.

\begin{thm}
Let $n\geq 1$ be an integer. There is a one-to-one correspondence
between s-extremal families $\mathcal{F}\subseteq 2^{[n]}$ of
Vapnik-Chervonenkis dimension $1$ with $supp(\mathcal{F})=[n]$,
$\cap_{F\in \mathcal{F}} F=\emptyset$ and directed edge-labelled
trees on $n+1$ vertices, all edges having a different label from
$[n]$. $\blacksquare$
\end{thm}

As a corollary one can prove the following statement.

\begin{cor}
There are $2^n (n+1)^{n-2}$ different s-extremal families
$\mathcal{F}\subseteq 2^{[n]}$ of Vapnik-Chervonenkis dimension at
most $1$ with $supp(\mathcal{F})=[n]$ and $\cap_{F\in \mathcal{F}}
F=\emptyset$.
\end{cor}
\textbf{Proof:} There are $(n+1)^{n-2}$ different edge labelled
undirected trees on $n+1$ vertices (see \cite{Pihurko}), all edges
having a different label from $[n]$ and each of these trees can be
directed in $2^n$ ways.
$\blacksquare$\\[0.1 cm]

Simple examples of s-extremal set systems are down-sets, i.e. set
systems $\mathcal{F}$ such that for all $i\in [n]$ $i\in F\in
\mathcal{F}$ implies $F\backslash\{i\}\in \mathcal{F}$. For
down-sets $Sh(\mathcal{F})=\mathcal{F}$, so they are obviously
s-extremal. One can obtain other examples from down-sets using
different set system operations, e.g. bit flips. For $i\in [n]$
let $\varphi_i$ be the the $i$th bit flip, i.e. for $F\in 2^{[n]}$

\[\varphi_i(F)=\left\{
\begin{array}{ll}
  \mbox{$F\backslash \{i\}$} & \mbox{ if $i\in F$} \\
  \mbox{$F\cup \{i\}$} & \mbox{ if $i\notin F$}.\\
\end{array}
\right.\]

For $\mathcal{F}\subseteq 2^{[n]}$ set
$\varphi_i(\mathcal{F})=\{\varphi_i(F)\ |\ F\in \mathcal{F}\}$. It
is easily seen that s-extremality is invariant with respect to
this operation since it keeps the family of shattered sets.
However not all s-extremal set systems can be obtained in this
way. For this note that in terms of the inclusion graph a bit flip
in the $i$th coordinate corresponds just to reversing the
direction of the edge with label $i$ in $G_{\mathcal{F}}$, i.e.
bit flips preserve the undirected structure of the inclusion
graph. Using this we can obtain many s-extremal examples not
coming from down-sets using bit flips. It is enough to pick a tree
that is not a star and consider a set system corresponding to any
possible orientation.

\section{s-extremal set systems from two consecutive layers}\label{2 layers}

\tab  For an uniform family $\mathcal{F}$ the graph
$G_{\mathcal{F}}$ is not connected, hence $\mathcal{F}$ cannot be
s-extremal. As a relaxation of uniformity we consider families
which belong to two consecutive layers of $2^{[n]}$. The next
proposition shows that extremal families among them are actually
special cases of the previously studied one.

\begin{prop}\label{VC 2 szint}
Let $\mathcal{F}\subseteq {[n] \choose k }\cup {[n] \choose k-1
}$, $n\geq k\geq 1$ be an s-extremal family of subsets of $[n]$.
Then we have $VC-dim (\mathcal{F})\leq 1$.
\end{prop}
\textbf{Proof}: For $n=2$ the statement can be verified by an easy
case analysis. For $n>2$ we can do induction on $k$. For $k=1$ the
statement is just trivial. Now suppose that $k>1$ and the result
holds for all values smaller than $k$. We prove that such an
s-extremal family cannot shatter a subset of size $2$. Suppose the
contrary, namely that $\mathcal{F}$ shatters for example
$\{1,2\}$. Let

\[\mathcal{F}_0^{(n)}=\{F\ |\ F\in \mathcal{F}\mbox{ and }n\notin F\}\subseteq {[n-1]\choose k}\cup
{[n-1]\choose k-1}\] and
\[\mathcal{F}_1^{(n)}=\{F\backslash\{n\}\ |\ F\in \mathcal{F}\mbox{ and }n\in F\}\subseteq {[n-1]\choose k-1}\cup
{[n-1]\choose k-2}.\]

Since $\mathcal{F}$ is s-extremal both $\mathcal{F}_0^{(n)}$ and
$\mathcal{F}_1^{(n)}$ must be s-extremal (it follows easily from
the proof of the Sauer inequality, see e.g. \cite{Shattering
news}) and for the shattered sets we have that

\[Sh(\mathcal{F})=Sh(\mathcal{F}_0^{(n)})\cup Sh(\mathcal{F}_1^{(n)})\cup \{F\cup\{n\}\ |\ F\in Sh(\mathcal{F}_0^{(n)})\cap Sh(\mathcal{F}_1^{(n)})\}.\]

Since $n>2$, by the induction hypothesis $\{1,2\}\in
Sh(\mathcal{F}_1^{(n)})$ cannot hold, thus we have $\{1,2\}\in
Sh(\mathcal{F}_0^{(n)})$. In this way we constructed an s-extremal
family with the same properties but on a smaller ground set.
Continuing this we get to an s-extremal family
$\mathcal{F}\subseteq {[k]\choose k}\cup {[k]\choose k-1}$ that
shatters $\{1,2\}$. However this is easily seen to be impossible,
because for any $F\in \mathcal{F}$ we have $|F\cap \{1,2\}|\geq
1$. This finishes
the proof. $\blacksquare$\\[0.1 cm]


Using essentially the same argument one can prove the following:

\begin{prop}\label{VC t szint}
Let $\mathcal{F}\subseteq {[n] \choose k }\cup {[n] \choose k-1
}\cup\dots \cup {[n] \choose k-t+1 }$, $n\geq k\geq t-1 \geq 1$ be
an s-extremal family of subsets of $[n]$. Then we have $VC-dim
(\mathcal{F})\leq t-1$. $\blacksquare$
\end{prop}

We return now to the situation when $\mathcal{F}\subseteq {[n]
\choose k }\cup {[n] \choose k-1 }$ for some $n\geq k\geq 1$ and
$supp(\mathcal{F})=[n]$, $\cap_{F\in \mathcal{F}} F=\emptyset$.
Proposition \ref{extr char VC-dim=1} says in this case that
$\mathcal{F}$ is s-extremal iff $G_{\mathcal{F}}$ (the undirected
version) is a tree and all labels on the edges are different. As
before, we also have that this tree has $n+1$ vertices and $n$
edges. Permuting the labels on the edges corresponds just to a
permutation of the ground set, so if we want to characterize
s-extremal set systems up to isomorphism we can freely omit the
labels from the edges.

Now suppose that we are given a tree $T$ on $n+1$ vertices having
$n$ edges. $T$ can also be viewed as a bipartite graph (since it
is acyclic, and so contains no odd cycles) with partition classes
$\mathcal{A},\mathcal{B}$. Direct all edges from $\mathcal{A}$ to
$\mathcal{B}$, and let $\mathcal{F}$ be as before the set system
this directed tree just defines. It is easily seen that we have
$\mathcal{F}\subseteq {[n] \choose k }\cup {[n] \choose k-1 }$,
where $k=|\mathcal{A}|$ and using the characterization of
s-extremal families we also get that $\mathcal{F}$ is s-extremal.
If we swap the role of $\mathcal{A}$ and $\mathcal{B}$ we get the
"dual" set system
\[\mathcal{F}'=\{[n]\backslash F\ |\ F\in \mathcal{F}\}\subseteq {[n] \choose n-k+1 }\cup {[n] \choose n-k },\] which is clearly also s-extremal using the same reasoning.

Summarizing the preceding discussion, we have the following:

\begin{thm}
Up to isomorphism and the operation of taking the "dual" of a set
system, there is a one to one correspondence between s-extremal
set systems $\mathcal{F}$ from two consecutive layers on the
ground set $[n]$ ($supp(\mathcal{F})=[n]$ and $\cap_{F\in
\mathcal{F}} F=\emptyset$) and trees on $n+1$ vertices. The
bijection is realized via the map $\mathcal{F}\rightarrow
G_{\mathcal{F}}$. $\blacksquare$
\end{thm}

\section{Ideal bases of s-extremal set systems of
$VC$-dimension at most $1$}\label{GR bases 1 VC dim}

Take a family $\mathcal{F}\subseteq 2^{[n]}$ and let $\F$ be a
field. For a set $F\subseteq [n]$ let $v_F\in \{0,1\}^n$ be its
characteristic vector, i.e. the the $i$th coordinate of $v_F$ is
$1$ exactly when $i\in F$. One can associate to $\mathcal{F}$ a
polynomial ideal $I(\mathcal{F})\lhd \Fx$, the vanishing ideal of
the set of characteristic vectors of the elements of
$\mathcal{F}$:

\[I(\mathcal{F})=\{f(x_1,\dots,x_n)\in \Fx\ |\ f(v_F)=0\mbox{ for all }F\in \mathcal{F}\}.\]

$I(\mathcal{F})$ carries a lot of information about the set
system. For this connection among $\mathcal{F}$ and
$I(\mathcal{F})$ see \cite{Springer} and \cite{Shattering news}.

If one is working with polynomial ideals it is advantageous to
have a good ideal basis. Now we briefly introduce one such class
of bases, namely Gr\"obner bases. For details we refer to
\cite{B1}, \cite{B2}, \cite{B3}, \cite{Cox}, and \cite{AL}.

A total order $\prec$ on the monomials composed from variables
$x_1,x_2,\dots, x_n$ is a term order, if $1$ is the minimal
element of $\prec$, and $\prec$ is compatible with multiplication
with monomials. One important term order is the lexicographic (lex
for short) order. We have $x_1^{w_1}\dots x_n^{w_n}\prec_{\rm
lex}x_1^{u_1}\dots x_n^{u_n}$ if and only if $w_i<u_i$ holds for
the smallest index $i$ such that $w_i\neq u_i$. Reordering the
variables gives an other lex term order.

The leading monomial $lm(f)$ of a nonzero polynomial $f\in\Fx$ is
the largest monomial (with respect to $\prec$) which appears with
nonzero coefficient in $f$, when written as the usual linear
combination of monomials. We denote the set of all leading
monomials of polynomials of a given ideal $I\lhd\Fx$ by $ Lm(I)
=\{ lm(f)\,:\,\,f\in I\}, $ and we simply call them the leading
monomials of $I$. A monomial is called a standard monomial of $I$,
if it is not a leading monomial of any $f\in I$. Let $Sm(I)$
denote the set of standard monomials of $I$. Obviously, a divisor
of a standard monomial is again in $Sm(I)$. Standard monomials
have some very nice properties, among other things they form a
linear basis of the $\F$-vector space $\Fx/I$.

A finite subset $\mathbb{G}\subseteq I$ is a Gr\"obner basis of
$I$, if for every $f\in I$ there exists a $g\in \mathbb{G}$ such
that $lm(g)$ divides $lm(f)$. It is not hard to verify  that
$\mathbb{G}$ is actually a basis of $I$, that is, $\mathbb{G}$
generates $I$ as an ideal of $\Fx$. It is a fundamental fact that
every nonzero ideal $I$ of $\Fx$ has a Gr\"obner basis (\cite{AL}
Corollary $1.6.5$).

Gr\"obner bases and standard monomials turned out to be very
useful when studying s-extremal set systems. Let
$\mathcal{F}\subseteq 2^{[n]}$ be a set system and
$I(\mathcal{F})$ its vanishing ideal. Subsets of $[n]$ can be
identified with square-free monomials via the map $H\rightarrow
\prod_{i\in H} x_i$. With this identification in mind one can
prove that $Sh(\mathcal{F})$, viewed as a set of monomials, is
just the union of the sets of standard monomials of
$I(\mathcal{F})$ for all lexicographic term orders.

\[Sh(\mathcal{F})=\bigcup_{lex\ term\ orders} Sm(I(\mathcal{F}))\]

On the other hand for vanishing ideals we have that
$|Sm(I(\mathcal{F}))|=|\mathcal{F}|$ for all term orders. These
facts altogether result that a set system $\mathcal{F}$ is
s-extremal iff the set of standard monomials is the same for all
lexicographic term orders. This algebraic characterization of
s-extremal set systems leads to an efficient algorithm for testing
s-extremality of a set system and offers also the possibility to
generalize the notion to arbitrary sets of vectors.  (For more
details and proofs see \cite{Springer}).

As an application of Proposition \ref{extr char VC-dim=1}, we
determine the Gr\"obner bases of s-extremal set systems of
$VC$-dimension $1$. Suppose that we are given a family
$\mathcal{F}\subseteq 2^{[n]}$ together with $Sh(\mathcal{F})$.
According to Subsection 4.2. of \cite{Springer} one can construct
a Gr\"obner basis of $I(\mathcal{F})$ as follows. For a pair of
sets $H\subseteq S\subseteq [n]$ define the following polynomial
\[
f_{S,H}=(\prod_{j\in H}x_j)(\prod_{i \in S\backslash H}(x_i-1)).
\]
Now if $S\notin Sh(\mathcal{F})$, then there exists a set $H
\subseteq S$ such that there is no set $F\in \mathcal{F}$ with
$F\cap S=H$. For this set $H$ we have $f_{S,H} \in
I(\mathcal{F})$. If the set $S$ is minimal (i.e. all proper
subsets $S'$ of $S$ are in $Sh(\mathcal{F})$) and $\mathcal{F}$ is
s-extremal, then we also have uniqueness for the corresponding
$H$. Moreover in the s-extremal case the collection of all these
$f_{S,H}$ polynomials corresponding to minimal elements outside
$Sh(\mathcal{F})$ together with $\{x_i^2-x_i, i \in [n]\}$ form a
Gr\"obner basis of $I(\mathcal{F})$ for all term orders. Actually
more is true:

\begin{prop}(\cite{Springer})\label{Springer}
$\mathcal{F}\subseteq 2^{[n]}$ is s-extremal iff there are
polynomials of the form $f_{S,H}$, which together with
$\{x_i^2-x_i, i \in [n]\}$ form a Gr\"obner basis of
$I(\mathcal{F})$ for all term orders. $\blacksquare$
\end{prop}

If we restrict ourselves to s-extremal set systems of
$VC$-dimension $1$, things become very simple. Suppose that
$\mathcal{F}\subseteq 2^{[n]}$ is a set system such that
$VC-dim(\mathcal{F})=1$, $supp(\mathcal{F})=[n]$ and $\cap_{F\in
\mathcal{F}} F=\emptyset$. In this case $Sh(\mathcal{F})$ is the
collection of all sets of size at most $1$, so the minimal sets
outside $Sh(\mathcal{F})$ are exactly the sets of size $2$. Fix
one such set $S=\{\alpha,\beta\}$, and consider the $2$ edges in
the inclusion graph $G_{\mathcal{F}}$ labelled by $\alpha$ and
$\beta$. From Proposition \ref{extr char VC-dim=1} we know that
$G_{\mathcal{F}}$ is a tree. Consider the unique path connecting
the $2$ edges. There are $4$ possibilities:

\begin{itemize}
\item The edges are directed towards each other on this path:

  
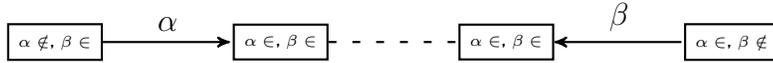
\begin{tikzpicture}[-,>=stealth',shorten >=1pt,auto,node distance=3cm,
    thick,main node/.style={rectangle,fill=none,draw,font=\sffamily\tiny\bfseries}]

    \node[main node] (1) {$\alpha \notin$, $\beta \in$};
    \node[main node] (2) [right of=1] {$\alpha \in$, $\beta \in$};
    \node[main node] (3) [right of=2] {$\alpha \in$, $\beta \in$};
    \node[main node] (4) [right of=3] {$\alpha \in$, $\beta \notin$};

    \draw[->] (1) -- (2) node[pos=.5,above] {$\alpha$};
    \draw[<-] (3) -- (4) node[pos=.5,above] {$\beta$};
    \draw[loosely dashed] (2) -- (3);
  \end{tikzpicture}

In this case the corresponding set $H$ is $\emptyset$, so
$f_{S,H}=(x_{\alpha}-1)(x_{\beta}-1)$. Indeed then every $F\in
\mathcal{F}$ contains either $\alpha$ or $\beta$.

\item The edges are directed away from each other each other on
this path:

   
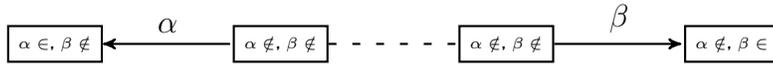
\begin{tikzpicture}[-,>=stealth',shorten >=1pt,auto,node distance=3cm,
     thick,main node/.style={rectangle,fill=none,draw,font=\sffamily\tiny\bfseries}]

     \node[main node] (1) {$\alpha \in$, $\beta \notin$};
     \node[main node] (2) [right of=1] {$\alpha \notin$, $\beta \notin$};
     \node[main node] (3) [right of=2] {$\alpha \notin$, $\beta \notin$};
     \node[main node] (4) [right of=3] {$\alpha \notin$, $\beta \in$};

     \draw[<-] (1) -- (2) node[pos=.5,above] {$\alpha$};
     \draw[->] (3) -- (4) node[pos=.5,above] {$\beta$};
     \draw[loosely dashed] (2) -- (3);
   \end{tikzpicture}

In this case the corresponding set $H$ is $\{\alpha,\beta\}$, so
$f_{S,H}=x_{\alpha}x_{\beta}$. No $F\in \mathcal{F}$ contains
$\{\alpha,\beta\}$.

\item The edges are directed in the same direction towards the
edge with label $\alpha$ on this path:

    
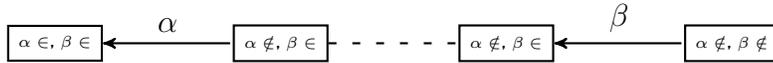
\begin{tikzpicture}[-,>=stealth',shorten >=1pt,auto,node distance=3cm,
      thick,main node/.style={rectangle,fill=none,draw,font=\sffamily\tiny\bfseries}]

      \node[main node] (1) {$\alpha \in$, $\beta \in$};
      \node[main node] (2) [right of=1] {$\alpha \notin$, $\beta \in$};
      \node[main node] (3) [right of=2] {$\alpha \notin$, $\beta \in$};
      \node[main node] (4) [right of=3] {$\alpha \notin$, $\beta \notin$};

      \draw[<-] (1) -- (2) node[pos=.5,above] {$\alpha$};
      \draw[<-] (3) -- (4) node[pos=.5,above] {$\beta$};
      \draw[loosely dashed] (2) -- (3);
    \end{tikzpicture}

In this case the corresponding set $H$ is $\{\alpha\}$, so
$f_{S,H}=x_{\alpha}(x_{\beta}-1)$. If $\alpha\in F$ for some $F\in
\mathcal{F}$ then $\beta\in F$ as well.

\item The edges are directed in the same direction towards the
edge with label $\beta$ on this path:

    
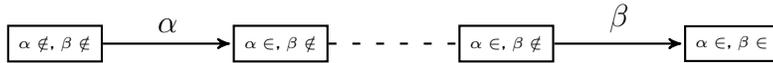
\begin{tikzpicture}[-,>=stealth',shorten >=1pt,auto,node distance=3cm,
      thick,main node/.style={rectangle,fill=none,draw,font=\sffamily\tiny\bfseries}]

      \node[main node] (1) {$\alpha \notin$, $\beta \notin$};
      \node[main node] (2) [right of=1] {$\alpha \in$, $\beta \notin$};
      \node[main node] (3) [right of=2] {$\alpha \in$, $\beta \notin$};
      \node[main node] (4) [right of=3] {$\alpha \in$, $\beta \in$};

      \draw[->] (1) -- (2) node[pos=.5,above] {$\alpha$};
      \draw[->] (3) -- (4) node[pos=.5,above] {$\beta$};
      \draw[loosely dashed] (2) -- (3);
    \end{tikzpicture}

Similarly to the previous case $H=\{\beta\}$, so
$f_{S,H}=(x_{\alpha}-1)x_{\beta}$.

\end{itemize}

Now if we have $G_{\mathcal{F}}$, then using the above case
analysis, one can easily compute a Gr\"obner basis for
$I(\mathcal{F})$. If we want just a basis of $I(\mathcal{F})$ and
not necessarily a Gr\"obner basis, we do not need to consider all
pairs. Consider $3$ consecutive edges in $G_{\mathcal{F}}$, i.e. a
path of length $3$ with labels $\alpha, \beta,\gamma$.

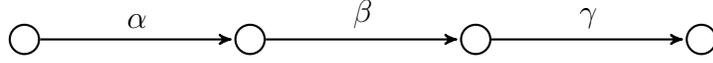
\begin{center}

    \begin{tikzpicture}[-,>=stealth',shorten >=1pt,auto,node distance=3cm,
      thick,main node/.style={circle,fill=none,draw,font=\sffamily\tiny\bfseries}]

      \node[main node] (1) {$$};
      \node[main node] (2) [right of=1] {$$};
      \node[main node] (3) [right of=2] {$$};
      \node[main node] (4) [right of=3] {$$};

      \draw[->] (1) -- (2) node[pos=.5,above] {$\alpha$};
      \draw[->] (2) -- (3) node[pos=.5,above] {$\beta$};
      \draw[->] (3) -- (4) node[pos=.5,above] {$\gamma$};

    \end{tikzpicture}

\end{center}

They define $3$ pairs and hence $3$ polynomials,
$f_{\alpha,\beta}=(x_{\alpha}-\varepsilon_{\alpha})(x_{\beta}-\varepsilon_{\beta})$,
$f_{\alpha,\gamma}=(x_{\alpha}-\varepsilon_{\alpha})(x_{\gamma}-\varepsilon_{\gamma})$,
$f_{\beta,\gamma}=(x_{\beta}-1+\varepsilon_{\beta})(x_{\gamma}-\varepsilon_{\gamma})$,
where $\varepsilon_{\alpha}$, $\varepsilon_{\beta}$ and
$\varepsilon_{\gamma}$ are $0$ or $1$ depending on the
orientations of the edges. However
\[(x_{\gamma}-\varepsilon_{\gamma})f_{\alpha,\beta}-(x_{\alpha}-\varepsilon_{\alpha})f_{\beta,\gamma}=(1-2\varepsilon_{\beta})f_{\alpha,\gamma},\]
where $1-2\varepsilon_{\beta}$ is either $1$ or $-1$, so
$f_{\alpha,\gamma}$ is superfluous, since it can be obtained from
$f_{\alpha,\beta}$ and $f_{\beta,\gamma}$. This means that when
constructing a basis of $I(\mathcal{F})$ it is enough to consider
only adjacent pairs of edges in
$G_{\mathcal{F}}$.

\section{s-extremal set systems of bounded $VC$-dimension}\label{t VC dim}

\tab The ideas from the previous sections can also be used to step
a bit further, and obtain results for s-extremal set systems of
bounded $VC$-dimension in general.\\

Let the projection of a set system $\mathcal{F}\subseteq 2^{[n]}$
to a set $X\subseteq [n]$ be

\[\mathcal{F}|_X=\{F\cap X\ |\ F\in \mathcal{F}\}.\]

Note that $X\in Sh(\mathcal{F})$ iff $\mathcal{F}|_X=2^X$.\\

The main result of this section considers the projections of a set
family from the perspective of extremality.

\begin{thm}\label{vetulet tetel}
Let $\mathcal{F}\subseteq 2^{[n]}$ be family of $VC$-dimension
$t\geq 1$ such that $\mathcal{F}|_X$ is s-extremal for all $2t+1$
element subsets $X$ of $[n]$. Let
\[\mathcal{G}=\{H\subseteq [n]\ |\ H\cap X\in \mathcal{F}|_X \mbox{ for all sets } X\subseteq [n] \mbox{ of size }2t+1\}.\]

Then $\mathcal{G}$ contains $\mathcal{F}$, $\mathcal{G}$ is
s-extremal and of $VC$-dimension $t$. Moreover, if $\mathcal{F}$
is s-extremal then we have $\mathcal{G}=\mathcal{F}$.
\end{thm}

Before proving Theorem \ref{vetulet tetel} we first present some
observations about extremal set systems in general.

For a set system $\mathcal{F}\subseteq 2^{[n]}$ and $i\in [n]$ let
$\mathcal{F}_0^{(i)}$ and $\mathcal{F}_1^{(i)}$ be as defined
previously in Section \ref{2 layers}. The downshift operation of a
set family $\mathcal{F}$ by the element $i\in [n]$ is defined as
follows:

\begin{center}
$D_i(\mathcal{F})=\{F\ |\ F\in \mathcal{F},\ i\notin F\}\cup\{F\
|\ F\in \mathcal{F},\ i\in F,\ F\backslash\{i\}\in
\mathcal{F}\}$\\
$\cup\{F\backslash\{i\}\ |\ F\in \mathcal{F},\ i\in F,\
F\backslash\{i\}\notin \mathcal{F}\}$\\
$=\{F\backslash\{i\}\ |\ F\in \mathcal{F}\}\cup\{F\ |\ F\in
\mathcal{F},\ i\in F,\ F\backslash\{i\}\in \mathcal{F}\}.$
\end{center}

It is not hard to see that $D_i$ preserves s-extremality (e.g.
\cite{defect sauer}, Lemma $1$) and as already noted above, if
$\mathcal{F}$ is s-extremal, then so is $\mathcal{F}_j^{(i)}$ for
$i\in [n]$ and $j=0,1$.

For a set system $\mathcal{F}\subseteq 2^{[n]}$ and $i\in [n]$ we
put

\begin{center}
$M_{i}(\mathcal{F})=\mathcal{F}_0^{(i)}\cap\mathcal{F}_1^{(i)}$,\\
$U_{i}(\mathcal{F})=\mathcal{F}_0^{(i)}\cup\mathcal{F}_1^{(i)}$.
\end{center}

\begin{prop}
Suppose that we are given a set system $\mathcal{F}\subseteq
2^{[n]}$ and an arbitrary element $i\in [n]$. Then if
$\mathcal{F}$ is s-extremal then so are $M_{i}(\mathcal{F})$,
$U_{i}(\mathcal{F})$ and $\mathcal{F}|_X$ for all $X\subseteq
[n]$. Actually for $\mathcal{F}|_X$ we have that
$Sh(\mathcal{F}|_X)=Sh(\mathcal{F})|_X$, more precisely a set
$Y\subseteq X$ is in $Sh(\mathcal{F}|_X)$ iff it is in
$Sh(\mathcal{F})$.
\end{prop}

\textbf{Proof:} The following equalities follow easily from the
definitions:

\begin{center}
$M_{i}(\mathcal{F})=\mathcal{F}_0^{(i)}\cap\mathcal{F}_1^{(i)}=D_{i}(\mathcal{F})_1^{(i)}$,\\
$U_{i}(\mathcal{F})=\mathcal{F}_0^{(i)}\cup\mathcal{F}_1^{(i)}=D_{i}(\mathcal{F})_0^{(i)}$.
\end{center}

From these it follows that if $\mathcal{F}$ is s-extremal then so
are $M_{i}(\mathcal{F})$ and $U_{i}(\mathcal{F})$, since we can
obtain them from $\mathcal{F}$ using operations preserving
s-extremality.

Next note that if $X=\{x_1,\dots,x_m\}$ then $\mathcal{F}|_X$ is
just $U_{x_1}(U_{x_2}(\dots U_{x_m}(\mathcal{F})\dots))$, thus if
the original set system is s-extremal, then so is its projected
version.

For the second part we only have to note that for some $Y\subseteq
X$ we have that $Y\cap (F\cap X)=Y\cap F$ for all $F\subseteq
[n]$, and the result follows.
$\blacksquare$ \\

We say that a set $I\subseteq [n]$ is strongly traced or strongly
shattered (\cite{reverse kleitman}, \cite{defect sauer}) by a set
system $\mathcal{F}\subseteq 2^{[n]}$ when there is a set
$B\subseteq [n]\backslash I$ such that

\[B+2^I=\{B\cup H\ |\ H\subseteq I\}\subseteq \mathcal{F}.\]

The collection of all sets strongly traced by $\mathcal{F}$ is
denoted by $st(\mathcal{F})$. It can be shown that
$|st(\mathcal{F})|$ is bounded from above by $|\mathcal{F}|$
(reverse Sauer inequality, \cite{reverse kleitman}, \cite{defect
sauer}), and a set system is called extremal with respect to the
reverse Sauer inequality if $|st(\mathcal{F})|=|\mathcal{F}|$. The
authors in \cite{defect sauer} proved that a set system is
extremal with respect to the original Sauer inequality exactly
when it is extremal with respect to the reverse one, and thus in
this case $Sh(\mathcal{F})=st(\mathcal{F})$.

Similarly to the case of $Sh(\mathcal{F})$, $st(\mathcal{F})$ can
also be obtained from the standard monomials of the vanishing
ideal $I(\mathcal{F})$, namely, if viewed as a set of monomials,
then it is the collection of those monomials which are standard
monomials for all lexicographic term orders.

For a set system $\mathcal{F}\subseteq 2^{[n]}$ and a set
$B\subseteq [n]$ denote the set family $\{I\subseteq [n]\backslash
B\ |\ I+2^B\subseteq \mathcal{F}\}$ by $\mathcal{F}(B)$. We remark
that if $B=\{i_1,\dots,i_m\}\subseteq [n]$, then $\mathcal{F}(B)$
is just $M_{i_1}(M_{i_2}(\dots M_{i_m}(\mathcal{F})\dots))$, and
hence is s-extremal if $\mathcal{F}$ is such. This together with
the fact that s-extremality of a set system implies the
connectedness of its inclusion graph proves in a simple way the
'if' direction of the following remarkable result of Bollobás and
Radcliffe (Theorem $3$ in \cite{defect sauer}).

\begin{thm}\label{extr <=> S(B) conn.}(\cite{defect sauer})
$\mathcal{F}\subseteq 2^{[n]}$ is s-extremal iff
$G_{\mathcal{F}(B)}$ is connected for every $B\subseteq [n]$.
$\blacksquare$
\end{thm}

%

Now we prepare the ground for the proof of Theorem \ref{vetulet
tetel} and return to the algebraic point of view. Let
$\mathcal{F}\subseteq 2^{[n]}$ be an arbitrary family, and fix one
term order $\prec$ on the monomials in $\Fx$. Suppose that we have
at our disposal a Gr\"obner basis $\mathbb{G}$ of $I(\mathcal{F})$
with respect to a term order $\prec$ (e.g. in the s-extremal case
we can compute one as described in Section \ref{GR bases 1 VC
dim}). From $\mathbb{G}$ we can compute a Gr\"obner basis for
$I(\mathcal{F}|_X)$: we only have to take the polynomials in
$\mathbb{G}$ depending only on the variables $x_i,\ i\in X$. In
general for a finite set of polynomials $\mathbb{G}\subseteq \Fx$
denote $\mathbb{G}\cap \mathbb{F}[x_i\ |\ i\in X]$ by
$\mathbb{G}|_X$. The leading term $lt(f)$ of a nonzero polynomial
$f\in \Fx$ with respect to $\prec$ is $lm(f)$ together with its
coefficient from $\F$. The S-polynomial of nonzero polynomials
$f,g$ in $\Fx$ is

\[S(f,g)=\frac{L}{lt(f)}f-\frac{L}{lt(g)}g,\]

\noindent where $L$ is the least common multiple of $lm(f)$ and
$lm(g)$. Buchberger's theorem (\cite{AL} Theorem $1.7.4.$) states
that a finite set $\mathbb{G}$ of polynomials in $\Fx$ is a
Gr\"obner basis for the ideal generated by $\mathbb{G}$ iff the
S-polynomial of any two polynomials from $\mathbb{G}$ can be
reduced to $0$ using $\mathbb{G}$. (For more
details on reduction and proofs see Chapter $1$ of \cite{AL}.)\\

\textbf{Proof of Theorem \ref{vetulet tetel}:} The facts
$\mathcal{F}\subseteq \mathcal{G}$ and $VC-dim(\mathcal{G})=t$
just follow immediately from the definition of $\mathcal{G}$.

Now fix one term order $\prec$ and one $2t+1$ element subset $X$
of $[n]$. Since by assumption $\mathcal{F}|_X$ is extremal,
according to Proposition \ref{Springer} the polynomials of the
form $f_{S,H}$, where $S$ is a minimal element outside
$Sh(\mathcal{F}|_X)$ and $H$ is the unique subset of $S$ such that
$H\notin (\mathcal{F}|_X)|_S$, together with $\{x_i^2-x_i\ |\ i\in
X\}$ form a Gr\"obner basis $\mathbb{G}_X$ of $I(\mathcal{F}|_X)$
with respect to $\prec$. We have

\[VC-dim(\mathcal{F}|_X)\leq VC-dim(\mathcal{F})=t,\]

\noindent hence $|S|\leq t+1$ for all polynomials of the form
$f_{S,H}$ in $\mathbb{G}_X$. Write

\[\mathbb{G}:=\bigcup_{X\subseteq [n], |X|=2t+1} \mathbb{G}_X.\]

Note that a polynomial $f_{S,H}$ from $\mathbb{G}$ is a member of
$\mathbb{G}_X$ for all $2t+1$ element subsets $X$ for which
$S\subseteq X$. First we prove using Buchberger's theorem that
$\mathbb{G}$ is a Gr\"obner basis of $\langle \mathbb{G}\rangle$
with respect to $\prec$. Take two polynomials $f,g\in \mathbb{G}$
and let $m$ be the number of variables occurring in them. If
$m\leq 2t+1$ then there is some $2t+1$ element set $X$ such that
$f,g\in \mathbb{G}_X$. However since $\mathbb{G}_X$ is a Gr\"obner
basis of $I(\mathcal{F}|_X)$, $S(f,g)$ can be reduced to $0$ using
$\mathbb{G}_X$, and so using $\mathbb{G}$ with respect to $\prec$
as well. On the other hand $m>2t+1$ is possible only if
$f=f_{S,H}$ and $g=f_{S',H'}$ for some appropriate sets
$H\subseteq S, H'\subseteq S'$ such that $S\cap S'=\emptyset$ and
$|S|=|S'|=t+1$. The leading terms of $f_{S,H}$ and $f_{S',H'}$ are
$\prod_{i\in S}x_i$ and $\prod_{j\in S'}x_j$ respectively, so we
can write them in the following form:

\[f=f_{S,H}=\underbrace{\prod_{i\in S}x_i}_{\ve x_S}+f'\mbox{ and }g=f_{S',H'}=\underbrace{\prod_{j\in S'}x_j}_{\ve x_{S'}}+g',\]

where $f'$ and $g'$ depend on disjoint sets of variables.

\[S(f,g)=S(f_{S,H},f_{S',H'})=\ve x_{S'}f'-\ve x_S g'\]

Here if we replace $\ve x_{S'}$ by $-g'$ and $\ve x_{S}$ by $-f'$,
the resulting polynomial will be identically $0$, i.e. reducing
$S(f,g)$ using $f,g\in \mathbb{G}$ gives $0$. Moreover, the above
reasoning works for all term orders $\prec$, so $\mathbb{G}$ is a
Gr\"obner basis of $\langle \mathbb{G}\rangle$ for all term
orders. Consider next

\[A:=\Fx / \langle x_i^2-x_i, i\in [n]\rangle.\]

Clearly we have $\{x_i^2-x_i\ |\ i\in[n]\}\subseteq \mathbb{G}$
and so $\langle\mathbb{G}\rangle$ defines also an ideal of $A$.
$A$ is actually isomorphic to the ring of all functions from
$\{0,1\}^n$ to $\F$, which in turn is isomorphic to $\F^{2^n}$. In
this ring every ideal is the intersection of maximal ideals and
hence every ideal is a radical ideal. This implies that every
ideal in $A$, in particular $\langle\mathbb{G}\rangle$ as well, is
a vanishing ideal of some finite point set from $\{0,1\}^n$. When
considering the $0-1$ vectors as characteristic vectors, this
finite point set also defines a finite set system. It is not
difficult to see, that in case of $\langle\mathbb{G}\rangle$ the
only possible candidate for this set system is $\mathcal{G}$
itself, so $\langle\mathbb{G}\rangle=I(\mathcal{G})$. However in
this case we get that $\mathbb{G}$ is a Gr\"obner basis of
$I(\mathcal{G})$ for all term orders and hence according to
Proposition \ref{Springer} $\mathcal{G}$ is s-extremal.

Finally we note that if $\mathcal{F}$ itself is already
s-extremal, then according to Proposition \ref{Springer}
$\mathbb{G}$ is a Gr\"obner basis for $I(\mathcal{F})$ as well and
so $\mathcal{F}=\mathcal{G}$.
$\blacksquare$\\[0.1 cm]

\section{Concluding remarks}

Concerning the structure of s-extremal set systems the question
arises whether an extremal family can be built up from the empty
system by adding sets to it  one-by-one in such a way that at each
step we have an s-extremal family.

\begin{conj}\label{bovit sejt}
For a nonempty s-extremal family $\mathcal{F}\subseteq 2^{[n]}$
does there always exist a set $F\in \mathcal{F}$ such that
$\mathcal{F}\backslash \{F\}$ is still s-extremal?
\end{conj}

From Theorem $2$ of \cite {defect sauer} we know that
$\mathcal{F}$ is s-extremal iff $2^{[n]}\backslash \mathcal{F}$ is
s-extremal, thus the above question has an equivalent form:

\begin{conj}
For an s-extremal family $\mathcal{F}\subsetneq 2^{[n]}$ does
there always exist a set $F\in \mathcal{F}$ such that
$\mathcal{F}\cup \{F\}$ is still s-extremal?
\end{conj}

There are several special cases when the answer appears to be
true:\\

1. If $\mathcal{F}$ is a nonempty down set ($F\in \mathcal{F}$ and
$H\subseteq F$ then $H\in \mathcal{F}$), then $\mathcal{F}$ is
extremal since $Sh(\mathcal{F})=\mathcal{F}$. Moreover in this
case if we omit a maximal element from $\mathcal{F}$ then it
remains still a down set  and so it will be still s-extremal.\\

2. If $\mathcal{F}$ is an extremal family of $VC$-dimension $1$,
then according to Proposition \ref{extr char VC-dim=1}, if we omit
a set corresponding to a leaf, i.e. to a vertex of degree $1$ in
$G_{\mathcal{F}}$, then the resulting set system will still be
extremal.\\

3. Anstee in \cite{anstee no triangle} constructed set systems
$\mathcal{F}\subseteq 2^{[n]}$, $|\mathcal{F}| = {n \choose 0}+{n
\choose 1}+{n \choose 2}$ without triangles, i.e. set systems with
the property, that for all $3$-element subsets $F$ we have that
$\mathcal{F}|_{F}$ does not contain all $2$-element subsets of
$F$. Note that in particular $VC-dim(\mathcal{F})$ is bounded from
above by $2$, hence we have that $Sh(\mathcal{F})\subseteq
{[n]\choose 0}\cup{[n] \choose 1}\cup{[n]\choose 2}$, implying
that $|Sh(\mathcal{F})|\leq {n \choose 0}+{n \choose 1}+{n \choose
2}$. Comparing the sizes of $\mathcal{F}$ and $Sh(\mathcal{F})$ we
obtain that such set systems are s-extremal.

Clearly any such set system $\mathcal{F}$ contains the extremal
subsystem ${[n] \choose 0}\cup{[n] \choose 1}$. For the remaining
part of these set systems Anstee's construction can be interpreted
in an inductive way as follows:
\begin{itemize}
\item $\mathcal{F}_1:={[n] \choose 0}\cup{[n] \choose 1}$ \item
For $k=2,3,\dots,n$ suppose we already constructed
$\mathcal{F}_{k-1}$. Let $\mathcal{G}_{k-1}$ be the collection of
all $k-1$-element sets in $\mathcal{F}_{k-1}$. Define $G_{k-1}$ to
be a graph, whose vertex set is $\mathcal{G}_{k-1}$ and there is
an edge between $A,B\in \mathcal{G}_{k-1}$ exactly when
$|A\vartriangle B|=2$. Take a spanning tree $T_{k-1}$ of
$G_{k-1}$.
\[\mathcal{F}_{k}:=\mathcal{F}_{k-1}\cup\{A\cup B\ |\ (A,B)\mbox{ is an edge of }T_{k-1}\}\]
\item $\mathcal{F}:=\mathcal{F}_n$
\end{itemize}

It is not hard to prove that when we add $A\cup B$, there will be
a unique new element that gets into the family of shattered sets,
namely $A\vartriangle B$, hence the resulting system after each
step will be s-extremal. Reversing it, if $\mathcal{F}$ is such an
example, then its elements can be deleted one-by-one in such a way
that the
remaining set system is s-extremal after each step.\\

4. More generally one can consider set systems
$\mathcal{F}\subseteq 2^{[n]}$ with the property, that for all
$t$-element subsets $F$ we have that $\mathcal{F}|_{F}$ does not
contain all $l$-element subsets of $F$, for some $l$ with $n\geq
t\geq l\geq 0$. Füredi and Quinn in \cite{Furedi} constructed for
all values $n\geq t\geq l\geq 0$ a set system $\mathcal{F}(n,t,l)$
with the desired property and of size
$|\mathcal{F}(n,t,l)|=\sum_{i=0}^{t-1} {n\choose i}$. The same
argument as above shows that $Sh(\mathcal{F}(n,t,l))$ consists of
all sets of size at most $t-1$ and hence $\mathcal{F}(n,t,l)$ is
s-extremal for all possible values. Their construction is as
follows.

For $x_1,\dots,x_i\in [n]$, $x_1<\dots<x_i$ let

\begin{center}
$E(x_1,\dots,x_i)=\{x\in[n]\ |\ x=x_j\mbox{ for } j\leq l \}$\\
$\cup\{x\in[n]\ |\ x>x_l\mbox{ but }x\neq x_j\mbox{ for any
}j>l\},$
\end{center}

\noindent in particular $E(\emptyset)=\emptyset$. Let
$\mathcal{F}(n,t,l)$ consist of all $E(x_1,\dots,x_i)$ where
$i\leq t-1$. Order the sets of $\mathcal{F}(n,t,l)$ as follows:
$E(X)\succ E(Y)$ if either $|X|>|Y|$, or $|X|=|Y|$ and $X\succ Y$
with respect to the standard lexicographic ordering. It is not
hard to see, that if we remove the elements of
$\mathcal{F}(n,t,l)$ with respect to this ordering one-by-one,
starting from the largest one, then each time when we remove some
$E(X)$, then $X$ is eliminated from the family of shattered sets,
hence after each step the resulting family will be still
s-extremal.\\

\begin{Large}\textbf{Acknowledgments}\end{Large}\\

We are grateful to Zolt\'an F\"uredi for discussions on the topic.


\end{document}